\documentclass[10pt]{amsart}
\RequirePackage{etex}
\usepackage[english]{babel}
\usepackage[utf8]{inputenc}
\usepackage{amsmath}
\usepackage{amssymb}
\usepackage{amsfonts}
\usepackage{amsthm}
\usepackage{mathrsfs}
\usepackage[all]{xy}
\usepackage[pdftex]{graphicx}
\usepackage{color}
\usepackage{cite}
\usepackage{url}
\usepackage{indent first}
\usepackage[labelfont=bf,labelsep=period,justification=raggedright]{caption}
\usepackage[english]{babel}
\usepackage[utf8]{inputenc}
\usepackage{hyperref}
\usepackage[capitalize,nameinlink]{cleveref}
\usepackage[colorinlistoftodos]{todonotes}
\usepackage{tkz-fct}
\usepackage[margin=1in]{geometry} 
\usepackage{tikz}
\usetikzlibrary{decorations.markings}

\DeclareMathOperator{\bip}{bip}

\newcommand{\NN}{\mathbb{N}}

\newcommand{\ind}{\mathbf{1}}
\newcommand{\ol}{\overline}

\theoremstyle{plain}
\newtheorem{thm}{Theorem}

\newtheorem{conj}[thm]{Conjecture}
\newtheorem{prop}[thm]{Proposition}

\newtheorem{obs}[thm]{Observation}

\theoremstyle{definition}
\newtheorem{defn}[thm]{Definition}

\theoremstyle{remark}
\newtheorem{rem}[thm]{Remark}

\numberwithin{equation}{section}
\numberwithin{thm}{section}

\begin{document}

\title{A note on directed analogues of the Sidorenko and forcing conjectures}

\author{Jacob Fox, Zoe Himwich, Nitya Mani, Yunkun Zhou}

\date{\today}

\begin{abstract} 
We study analogues of Sidorenko's conjecture and the forcing conjecture in oriented graphs, showing that natural variants of these conjectures in directed graphs are equivalent to the asymmetric, undirected analogues of the conjectures. 

\end{abstract}

\maketitle

\section{Introduction}
Estimating the minimum possible number of copies of a graph $A$ in another graph $H$ with a given number of vertices and edges is a central problem in extremal graph theory. 
One of the most important open problems in this area is \textit{Sidorenko's conjecture}. Informally, Sidorenko's conjecture states that every dense graph $G$ has asymptotically at least as many copies of any fixed bipartite graph as is expected in a random graph with the same number of vertices and the same edge density as $H$.
This conjecture was independently posed by Erd\H{o}s-Simonovits \cite{SIM82} and Sidorenko \cite{SID93}
and often appears in the language of \textit{graph homomorphisms}, vertex maps of graphs that send edges to edges. We use the standard notation for a graph 
$H$, that its vertex set is denoted $V(H)$, its edge set $E(H)$, its number of vertices is $v(H)=|V(H)|$, and its number of edges is $e(H)=|E(H)|$.

\begin{defn}\label{def:hom-density}
Given undirected graphs $A, H$, we call $f: V(A) \to V(H)$ a \textit{homomorphism} if for $x, y\in V(A)$, $f(x)$ is adjacent to $f(y)$ whenever $x$ is adjacent to $y$. Let $h(A, H)$ be the number of homomorphisms $V(A) \rightarrow V(H)$. The \textit{homomorphism density} of $A$ in $H$ is the fraction of vertex maps $A \rightarrow H$ that are homomorphisms, given by $t(A, H) := h(A, H)/v(H)^{v(A)}.$
\end{defn}


\begin{conj}[Sidorenko's Conjecture]\label{c:sidconj}
For every bipartite graph $A$ and graph $H$, we have $$t(A, H) \ge t(K_2, H)^{e(A)}.$$
\end{conj}

Sidorenko~\cite{SID93} showed that this conjecture holds for several types of graphs, including complete bipartite graphs, trees, and even cycles. While Sidorenko's conjecture remains open, a variety of special cases have been resolved. Conlon, Fox, and Sudakov~\cite{CON10} showed that the conjecture holds for bipartite graphs with one vertex complete to the other side. In several other works, including~\cite{LI17,KIM16,CON18,SZE}, further progress has been made to show the conjecture holds for larger subfamilies of the collection of bipartite graphs. The M\"obius ladder on 10 vertices, the undirected graph obtained by deleting the edges of a 10-cycle from a $K_{5,5},$ is a notorious small open case of Sidorenko's conjecture.

Sidorenko's conjecture has a wide variety of applications to the study of random matrix theory, Markov chains, and to the study of \textit{quasirandomness}. Quasirandom graphs were first studied by Thomason \cite{TH87} and Chung-Graham-Wilson \cite{FAN89}; they observed that a large number of properties that Erd\H{o}s-Renyi random graphs satisfy are actually equivalent. This motivated the characterization of deterministic graph families that satisfy these equivalent properties. For $p \in [0, 1]$, we say that a sequence of distinct undirected graphs $\{H_n\}_{n \in \NN}$ is \textit{$p$-quasirandom} if $t(K_2, H_n) \to p$, and for all fixed undirected graphs $A$, $t(A, H_n) \to p^{e(A)}$ as $n \to \infty$.




There is a strengthening of Sidorenko's conjecture on characterizing quasirandom graph properties. A graph $A$ is \textit{$p$-forcing} if for all families of graphs $\{H_n\}_{n = 1}^{\infty}$, the family $\{H_n\}_{n = 1}^{\infty}$ is $p$-quasirandom if and only if the density of $H_n$ is asymptotically $p$, and the number of copies of $A$ in $H_n$ is asymptotically the number expected in the Erd\H{o}s-Renyi graphs $G(v(H_n), p)$. We call $A$ \textit{forcing} if it is $p$-forcing for all $p \in [0, 1]$. 
The \textit{forcing conjecture}, initially posed by Skokan and Thoma \cite{SKO04} states that subgraphs are forcing if and only if they are bipartite and contain a cycle (showing these conditions are necessary is straightforward).

While there has been an extensive effort over the past decades to resolve parts of Sidorenko's conjecture and related problems, the analogues of these problems for oriented graphs remain relatively poorly understood. The first discussion of directed quasirandomness appeared in the context of tournaments in~\cite{FAN91}. Recently, substantial progress has been made in understanding oriented subgraph counts and quasirandomness in tournaments and more general oriented graphs as in~\cite{GRI13,COR15, COR17, BU19, HA19,FHMZ22}. 

In this article, we investigate natural analogues of the Sidorenko and forcing properties in oriented graphs. Throughout, we work with \textit{oriented graphs}, which are formed from (simple) undirected graphs by orienting each edge. We can similarly define homomorphism density for directed graphs as in~\cref{def:hom-density}, where $f: V(B)\to V(G)$ is a homomorphism if for $x, y\in V(B)$, we have $(x, y)\in E(B)$ if and only if $(f(x), f(y))\in E(G)$. An oriented graph $B$ has the \textit{directed Sidorenko property} if for any oriented graph $G$ with edge density $p$, the number of copies of $B$ in $G$ is at least as many as expected if we orient each edge of an Erd\H{o}s-Renyi random graph with $v(G)$ vertices and density $p$ uniformly at random. 

\begin{defn}[Directed Sidorenko]~\label{d:disidorenko}
An oriented graph $B$ is said to have the \textit{directed Sidorenko property} if for every oriented graph $G$,
\begin{equation}\label{e:sid1}
t(B, G) \ge t(\vec K_2, G)^{e(B)}, 
\end{equation}
where $\vec K_2$ refers to the oriented complete graph on two vertices.
\end{defn} 

For a bipartite oriented graph to have the directed Sidorenko property, it must have a homomorphism to an edge (see~\cref{t:disidorenkonec}). We conjecture that this necessary property is in fact sufficient. 

\begin{conj}\label{c:sidsid}
If $B$ is a bipartite oriented graph with a homomorphism $B\rightarrow \vec K_{2}$, then $B$ has the directed Sidorenko property.
\end{conj}

We prove that~\cref{c:sidsid} is equivalent to a previously studied conjecture. A bipartite undirected graph $A$ has the asymmetric Sidorenko property if $t_{\bip}(A, H)\geq p^{e(A)}$ for all bipartite undirected graphs $H=(U_{1}\sqcup U_{2},F)$ with bipartite density $p=\frac{|F|}{|U_{1}||U_{2}|}$ (\cref{d:asymsid}). Throughout this work, when we say $A$ is a bipartite graph, it also fixes a bipartition $V(A) = A_1\sqcup A_2$ of its vertices, where we refer to $A_1$ as the first part and $A_2$ the second part. Here, $t_{\bip}(A, H)$ is the fraction of vertex maps respecting the corresponding bipartitions, i.e. sending vertices in the $A_1$ to $U_1$ and $A_2$ to $U_2$, that are homomorphisms. 

Throughout this work, we let $\overline{B}$ denote the \textit{underlying undirected graph} of oriented graph $B$, the undirected graph obtained by forgetting the directions of all of the edges in $G$.
\begin{thm}\label{t:disidorenko}
If $B$ is an oriented graph with a homomorphism $B \rightarrow \vec K_2$, then $\overline{B}$ satisfies the asymmetric Sidorenko conjecture if and only if $B$ has the directed Sidorenko property. 
\end{thm}

We also exhibit two necessary conditions for an oriented graph to have the directed forcing property (defined precisely in~\cref{d:diforcing}).

\begin{thm}\label{t:forcingn} 
If an oriented graph $B = (V, E)$ has the directed forcing property, then it has a homomorphism to an $\vec K_2$ and $\overline{B}$ has a cycle.
\end{thm}

We conjecture the two necessary conditions above characterize oriented graphs with the directed forcing property. 
\begin{conj}\label{c:forc}
If $B$ is an oriented graph with a homomorphism to an edge and whose underlying undirected graph has a cycle, then $B$ has the directed forcing property.
\end{conj}

We also obtain an analogous result to our reduction from the asymmetric Sidorenko conjecture for oriented graphs with the directed forcing property.

\begin{thm}\label{t:diforcing}
Let $B$ be any oriented graph with a homomorphism $B \rightarrow \vec K_2$. $\overline{B}$ has the asymmetric forcing property if and only if $B$ has the directed forcing property. 
\end{thm} 

The asymmetric Sidorenko property is a stronger notion, as it implies the Sidorenko property. This can be seen by considering an auxiliary bipartite graph $H'$ constructed from some undirected graph $H$ by making two copies of the vertex set of $H$ and connecting pairs of vertices in $H'$ between the two parts if they correspond to edges in $H$. Then, for any undirected $A$ with the asymmetric Sidorenko property, by counting homomorphic copies of $A$ in $H'$ that send vertices to corresponding parts, we can find the desired Sidorenko lower bound on the number of homomorphisms of $A$ in $H$. 

It is unknown if the asymmetric and classical Sidorenko properties are equivalent (although this is conjectured to be the case). However, although there are no known black box reductions from the standard Sidorenko conjecture to the asymmetric Sidorenko conjecture, essentially all graphs known to have the Sidorenko property are also known to satisfy the asymmetric Sidorenko property (for example, in~\cite{CON10} it was observed that any bipartite graph with a vertex complete to one part has the asymmetric Sidorenko property).

We begin with some preliminaries on the directed and asymmetric Sidorenko properties along with graphons in~\cref{s:prelim}. Using the language of graphons, we prove~\cref{t:disidorenko} in~\cref{s:sid}. We subsequently study the directed forcing property in~\cref{s:force}. We conclude with some open problems and further directions in~\cref{s:conc}.





\section{Preliminaries}\label{s:prelim}
Throughout, we let oriented graph $G = (V, E)$ have vertex set $V$ of size $v(G)$ and edge set $E$ of size $e(G)$ with density $p = e(G)/\binom{v(G)}{2}$ and underlying undirected graph $\overline{G}$.
\begin{defn}
Given an oriented graph $G$, for vertex $v \in V$, its \textit{out-neighborhood} is denoted by $N^+(v) = \{w : (v, w) \in E\}$ and its \textit{out-degree} is $d^+(v) = |N^-(v)|$. We analogously define the \textit{in-neighborhood} $N^-(v)$ and associate to $v$ its \textit{in-degree} $d^-(v)$.
\end{defn}

\subsection{Directed and asymmetric Sidorenko}
Recall from~\cref{d:disidorenko} that oriented graph $B$ has the directed Sidorenko property if it is systematically over-represented in oriented graphs.



\begin{rem}
We can alternatively define $B$ to satisfy the directed Sidorenko property if 
for every oriented graph $G$ on $n$ vertices,
\begin{equation}\label{e:sid3}
N_L(B, G) \ge \left( \left( \frac{e(G)}{v(G)^2} \right)^{e(B)} - o(1) \right) v(G)^{v(B)},
\end{equation} where $N_L(B, G)$ denotes the number of labeled copies of $B$ in $G$ and $o(1)$ is a quantity that goes to $0$ as $v(G) \to \infty$.
It is not hard to show that this alternative definition is equivalent to that given by ~\cref{d:disidorenko}.

\end{rem}

\begin{rem}
If $G$ is a random oriented graph, where each pair of vertices appears as an edge with constant probability $p$ oriented uniformly at random (independently), an injective map $ V(B) \mapsto V(G)$ yields a copy of $B$ in $G$ with probability $(p/2)^{e(B)}$. Thus, Equation~\eqref{e:sid3} says that the number of copies of $B$ in $G$ is asymptotically at least a natural random bound. 
\end{rem}

We will reduce the directed Sidorenko property to an undirected, bipartite Sidorenko property.


\begin{defn}\label{d:asymsid}
Let $A=(V_{1}\sqcup V_{2}, E)$ and $H=(U_{1}\sqcup U_{2},F)$ be bipartite undirected graphs.
Let $t_{\bip}(A, H)$ be the density of maps $f:V(A)\to V(H)$ where $f(V_{i})\subset U_{i}$ for $i=1,2$ that are homomorphisms.
We say that $A$ has the \textit{asymmetric Sidorenko property} if $t_{\bip}(A, H) \ge t_{\bip}(K_2, H)^{e(A)}$ holds for all bipartite undirected graphs $H$. 
\end{defn}
Note that in the above definition, $H$ has bipartite edge density $t_{\bip}(K_2, H)=\frac{|F|}{|U_{1}||U_{2}|}$.

\subsection{Graphons}

The relationship between asymmetric Sidorenko/forcing properties and their oriented analogues is most naturally seen by leveraging the language of graphons, as in~\cite{LOV12}.  
The language of graphons will also give us a clean framework to understand quasirandomness and forcing in oriented and bipartite graphs. The study of quasirandomness in oriented graphs began with Chung and Graham~\cite{FAN91}, who studied 
quasirandom tournaments. Recently, several researchers have studied more general quasirandom oriented graphs~\cite{GRI13, AM11}. For conciseness, below we define quasirandomness and forcing \textit{only} in the setting of graphons; however, the notions we define below specify to the classical notions of quasirandomness in oriented and undirected bipartite graphs, and to forcing in bipartite undirected graphs.


\begin{defn}\label{def:asym-graphon}
An \textit{graphon} is a measurable function $W : [0, 1]^2 \rightarrow [0, 1]$.
Given a directed graph $G = ([n], E),$ we define its associated \textit{directed graphon} to be the function $W_G: [0, 1]^2 \rightarrow [0, 1]$ obtained by equipartitioning $[0, 1] = I_1 \sqcup \cdots \sqcup I_n$ into equal length intervals, and letting $W_G(x, y) = \ind\{(i, j) \in E\}$ for $x \in I_i, y \in I_j$ for $1\leq i, j\leq n$. 

Given a bipartite graph $H = (V \sqcup W, F)$ with $V = \{v_1, \ldots, v_n\}, W = \{w_1, \ldots, w_m\}$, we define its associated \textit{bipartite graphon} $W_H: [0, 1]^2 \rightarrow [0, 1]$ as following. We equipartition the interval $[0, 1]$ into equal length intervals $[0, 1] = I_1 \sqcup \cdots \sqcup I_n$ and $[0, 1] = J_1\sqcup \cdots \sqcup J_m$, and let $W_H(x, y) = \ind\{(v_i, w_j) \in F\}$ for $x \in I_i, y \in J_j$ for $1\leq i\le n, 1\le j \le m$.
\end{defn}
Symmetric graphons (which satisfy $W(x,y)=W(y,x)$) are commonly studied, as they correspond to undirected graph limits. Here, we study graphons more generally (in the \textit{asymmetric} setting) as limit objects of sequences of either directed graphs or bipartite undirected graphs with specified parts.

The space of graphons is compact under the cut metric, given by identifying graphons with cut distance zero, as defined below.

\begin{defn}
Given a measurable $W: [0, 1]^2 \rightarrow [0, 1],$ the \textit{cut norm} of $W$ is given by 
$$\|W\|_{\square} = \sup_{S \times T \subset [0, 1]^2} \left| \int_{S \times T} W(x, y) dx dy \right|,$$
taking the supremum over measurable $S, T.$
Given two graphons $W, U$, we define their \textit{cut distance} as
$$\delta_{\square}(U, W) = \inf_{\phi:[0,1] \rightarrow [0, 1]} \|U - W^{\phi}\|,$$
taking the infimum over measure preserving maps $\phi$, where $W^\phi(x, y) = W(\phi(x), \phi(y))$.
\end{defn}

In terms of cut distance, graphons are the limiting objects of graphs. Formally we have the following result (see e.g.~\cite{LOV12} for proofs and more detailed exposition about graphons).
\begin{thm}[Corollary 11.15~\cite{LOV12}]\label{thm:graphon-lim}
Suppose that $W:[0, 1]^2\to[0,1]$ is a measurable function. Then there exists a sequence of bipartite graphs $\{H_n\}_{n\in \NN}$ such that
\[\lim_{n\to\infty} \delta_\square(W, W_{H_n}) = 0.\]
If $W(x, y) + W(y, x) \le 1$ for all $x, y\in [0, 1]$, there is a sequence of oriented graphs $\{G_n\}_{n\in \NN}$ with
\[\lim_{n\to\infty} \delta_\square(W, W_{G_n}) = 0.\]
\end{thm}

The cut norm also offers a succinct characterization of quasirandomness for graphons.
\begin{defn}
A sequence of graphons $\{W_n\}_{n \in \NN}$ is \textit{$p$-quasirandom} if $$\lim_{n \rightarrow \infty} \|W_n - p\|_{\square} = 0.$$
\end{defn}

Quasirandomness of a family of oriented or bipartite graphs is equivalent to quasirandomness of the associated graphon family (c.f. \S1.4, Examples 11.37-38 in~\cite{LOV12}). 

The above definition of quasirandomness appears at first to have little to do with counts of oriented subgraphs. This relationship becomes more apparent via the below characterization (culminating in~\cref{lem:counting}) of graphon quasirandomness via homomorphism density of subgraphs. 


\begin{defn}\label{def:graphon-hom-density}
Given graphon $W$ and directed graph $B$, the \textit{$B$-density in $W$} is
$$t(B, W) = \int_{[0, 1]^{V(B)}} \prod_{(i, j) \in E(B)} W(x_i, x_j) \prod_{i \in V(B)} dx_i$$
Given graphon $W$ and bipartite undirected graph $A = (A_1 \sqcup A_2, E(A)),$ the \textit{bipartite $A$-density in $W$} is 
$$t_{\bip}(A, W) = \int_{[0, 1]^{A_1\sqcup A_2}} \prod_{(v_i, w_j) \in E(A)} W(x_i, y_j) \prod_{v_i \in A_1} dx_i \prod_{w_j \in A_2} dy_j.$$
A sequence of graphons is \textit{left convergent} to graphon $W$ if for every (bipartite/directed) graph $F$, $t(F, W_n) \rightarrow t(F, W)$ as $n \rightarrow \infty.$
\end{defn}

The above definitions have the following consequence.
\begin{obs} \label{lem:switching}
Let $\overline{B} = (V_1\sqcup V_2, E')$ be an undirected bipartite graph, and let $B = (V, E)$ be the directed graph obtained from $\overline{B}$ by directing all edges from $V_1$ to $V_2$. Then $t(B, W) = t_{\bip}(\overline{B}, W)$ for any graphon $W: [0, 1]^2\to [0, 1]$.
\end{obs}
\cref{lem:switching} follows by examining the expression for the homomorphism densities given in~\cref{def:graphon-hom-density} and noticing that the integrands are the same because the edges of $B = (V_1 \sqcup V_2, E(B))$ are exactly those of the form $(x, y)$ for $x \in V_1, y \in V_2$ and where $(x, y) \in E(\overline{B})$.

The above notions of homomorphism density are consistent with the analogous definitions in graphs in the following sense. For any directed graphs $B$ and $G$, $t(B, G) = t(B, W_G)$ where $W_G$ is the directed graphon associated to $G$; for any bipartite graphs $A$ and $H$, $t_{\bip}(A, H) = t_{\bip}(A, W_H)$ with $W_H$ the bipartite graphon associated to $H$.
\begin{prop}[\S11,~\cite{LOV12}]\label{prop:graphon-quasi}
A sequence of bipartite graphs $\{H_n\}_{n\geq 1}$ is $p$-quasirandom if and only if 
\[\lim_{n\to\infty} \|W_{H_n} - p\|_{\square} = 0.\]
A sequence of directed graphs $\{G_n\}_{n\geq 1}$ is $p$-quasirandom if and only if 
\[\lim_{n\to\infty} \|W_{G_n} - p\|_{\square} = 0.\]
\end{prop}

By identifying graphons that differ on a measure $0$ set, the space of graphons becomes a compact space with the cut distance metric (c.f. Theorem 9.23,~\cite{LOV12}).
As shown in~\cite{BCL08}, left convergence is equivalent to convergence in cut distance; a sequence of graphons $\{W_n\}_{n \in \NN}$ left converges to $W$ if and only if $\delta_{\square}(W_n, W) \rightarrow 0$ as $n \rightarrow \infty$. Quantitatively, we have the following counting lemma (which establishes one direction of the equivalence).

\begin{thm}[Lemma 10.23,~\cite{LOV12}]\label{lem:counting}
Let $U, W$ be graphons. If $B$ is an oriented graph, then $$|t(B, W) - t(B, U)| \le e(B) \delta_{\square}(U, W).$$
If $A$ is a bipartite undirected graph, we similarly have $$|t_{\bip}(A, W) - t_{\bip}(A, U)| \le e(A) \delta_{\square}(U, W).$$
\end{thm}

This equivalence establishes a homomorphism density version of quasirandomness as an immediate consequence; it also allows us to give a simple definition of the directed forcing property in the language of graphons.



\begin{defn}\label{d:diforcing}
Oriented graph $B$ has the \textit{directed forcing property} if for all $p \in [0, 1]$ and any sequence of oriented graphs, $\{G_n\}_{n \in \NN}$ with $t(\vec{K_2}, G_n) \to p/2$ and $t(B, W_{G_n}) \to \left(\frac{p}{2}\right)^{e(B)}$ as $n \to \infty$, then $\{G_n\}_{n \in \NN}$ is $p$-quasirandom.
Bipartite undirected graph $A$ has the \textit{asymmetric forcing property} if for all $p \in [0, 1]$ and any sequence of bipartite undirected graphs, $\{H_n\}_{n \in \NN}$ with bipartite density tending to $p$ as $n \to \infty$ and $t_{\bip}(A, W_{H_n}) \to p^{e(B)}$, $\{H_n\}_{n \in \NN}$ is $p$-quasirandom.
\end{defn}

\section{Directed Sidorenko property}\label{s:sid}
Which oriented graphs have the directed Sidorenko property? We reduce this question to the analogous conjecture in the asymmetric, undirected setting.

Let $G = K_{n,n}^{\rightarrow} = (V_1 \sqcup V_2, E)$ be a complete balanced bipartite oriented graph on $2n$ vertices where all edges $e = (v_1, v_2) \in E$ are oriented so that $v_1 \in V_1$ and $v_2 \in V_2$. Notice that if oriented graph $B$ does not have a homomorphism to an edge, there are zero copies of $B$ in $G$. Thus, in order to have the directed Sidorenko property, an oriented graph $B$ must have a homomorphism to an edge, giving the following.

\begin{obs}\label{t:disidorenkonec}
If oriented graph $B$ has the directed Sidorenko property, then it has a homomorphism $B \rightarrow \vec K_2$.
\end{obs}

Therefore, we only need to analyze directed graphs $B$ obtained from a bipartite graph $\ol B$ by orienting all edges from the first part to the second. We aim to show~\cref{t:disidorenko}, which says that $B$ is directed Sidorenko if and only if $\overline B$ is asymmetric Sidorenko. We first give the following characterization of the Sidorenko properties using graphons.

\begin{prop}\label{prop:diSidgraphon}
An oriented graph $B$ satisfies the directed Sidorenko property if and only if for all measurable functions $W: [0, 1]^2\to [0, 1]$, we have
\[t(B, W) \geq \left|\int_{[0, 1]^2}W(x, y)\,dxdy\right|^{e(B)} = t(\vec{K_2}, W)^{e(B)}.\]
\end{prop}
\begin{proof}
First, if $B$ satisfies $t(B, W) \geq t(\vec{K_2}, W)^{e(B)}$ for all $W$, then for any oriented graph $G$, we take $W = W_G$ and derive that $t(B, G) \geq (e(G)/v(G)^2)^{e(B)}$. This means that $B$ satisfies the directed Sidorenko property.

On the other hand, suppose that $B$ satisfies the directed Sidorenko property. As a consequence, for any $W: [0, 1]^2\to [0, 1]$, by~\cref{thm:graphon-lim}, there exists a sequence of directed graphs $\{G_n\}_{n\geq 1}$ such that $\lim_{n\to\infty}\delta_\square(W/2, W_{G_n}) = 0$. Since $B$ satisfies the directed Sidorenko property, we have
\[t(B, W_{G_n}) = t(F, G_n) \geq t(\vec{K_2}, W_{G_n})^{e(B)}.\]
Then by~\cref{lem:counting}, we have $t(B, W/2) \geq t(B, W_{G_n}) - e(B)\delta_\square(W/2, W_{G_n})$ and $t(\vec{K_2}, W_{G_n}) \geq t(\vec{K_2}, W/2) - \delta_\square(W/2, W_{G_n}).$ Hence we have
\[t(B, W/2) \geq \left|\int_{[0, 1]^2}\frac12 W(x, y)\,dxdy - \delta_\square(W/2, W_{G_n})\right|^{e(B)} - e(B)\delta_\square(W/2, W_{G_n}).\]
Since this is true for all $n$, as $n \to \infty$, we get $t(B, W/2) \geq t(\vec{K_2}, W/2)^{e(B)}$.
Multiply both sides by $2^{e(B)}$ we have the desired inequality, thereby establishing the equivalence.
\end{proof}
Using the same argument, we can prove the following analogue for bipartite graphs.
\begin{prop}\label{prop:bipSidgraphon}
An undirected bipartite graph $A$ satisfies the asymmetric Sidorenko property if and only if for all measurable functions $W: [0, 1]^2\to [0, 1]$, we have
\[t_{\bip}(A, W) \geq \left|\int_{[0, 1]^2}W(x, y)\,dxdy\right|^{e(A)}.\]
\end{prop}

Now we are ready to prove~\cref{t:disidorenko}.
\begin{proof}[Proof of~\cref{t:disidorenko}]
Let $\ol{B}$ be an undirected bipartite graph, and $B$ be the directed graph obtained from $\ol{B}$ by orienting all edges from the first part to the second. We aim to show that $\ol{B}$ is asymmetric Sidorenko if and only if $B$ is directed Sidorenko.
By~\cref{prop:bipSidgraphon}, $\ol{B}$ is asymmetric Sidorenko if and only if for all measurable $W:[0, 1]^2\to [0, 1]$ we have $$t_{\bip}(\ol B, W) \geq \left|\int_{[0, 1]^2}W(x, y)\,dxdy\right|^{e(\ol B)}.$$ By~\cref{lem:switching}, we have $t_{\bip}(\ol B, W) = t(B, W)$. As $e(\ol B) = e(B)$, we have $\ol B$ is asymmetric Sidorenko if and only if $$t(B, W) \geq \left|\int_{[0, 1]^2}W(x, y)\,dxdy\right|^{e(B)}$$ for all measurable $W:[0, 1]^2\to [0, 1]$. This is, by~\cref{prop:diSidgraphon}, equivalent to $B$ being directed Sidorenko. Hence we conclude that the two conditions are indeed equivalent.
\end{proof}

\section{Directed forcing property}\label{s:force}
We first show that, in order to have the directed forcing property, an oriented graph must have a homomorphism to $\vec K_2$.
\begin{prop}
If an oriented graph $B$ satisfies the directed forcing property, then it has a homomorphism $B\to\vec K_2$.
\end{prop}
\begin{proof}
We may assume that $B$ has no isolated vertices (by removing any isolated vertices that exist).
For $\lambda \in [0, 1]$, we define  graphon $W^{(\lambda)}$ as below and in~\cref{fig:wlambda}
\[W^{(\lambda)}(x, y) = \begin{cases}1 - \lambda & \textrm{if }x\in [1/4, 1/2], y\in [0, 1/4],\\
\lambda/4 & \textrm{if }x\in [1/2, 1], y\in [1/2, 1],\\
0 & \textrm{otherwise}.\end{cases}\]

\begin{figure}[ht!]
\begin{tikzpicture}
\draw[step=2cm,color=black] (-2.001,-2.001) grid (2,2);
\draw[step=1cm,color=black] (-2,0) grid (0,2);
\foreach \x/\y/\m in {-1.5/+1.5/$0$,-1.5/+0.5/$0$,-0.5/+0.5/$0$,-0.5/+1.5/\small{$1-\lambda$},+1/-1/$\frac{\lambda}{4}$,-1/-1/$0$,+1/+1/$0$} 
    \node at (\x,\y) {\m};
\node at (-2.25,2.25) {\tiny{$(0, 0)$}};
\node at (2.25,-2.25) {\tiny{$(1, 1)$}};
\end{tikzpicture}
\caption{Visual depiction of the graphon $W^{(\lambda)}$.}
\label{fig:wlambda}
\end{figure}
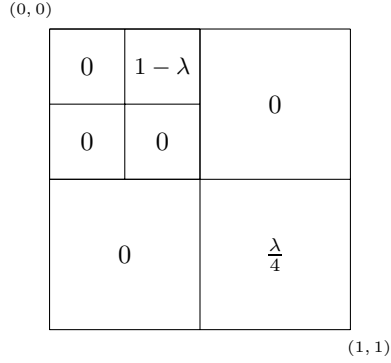

For any $\lambda$, we know that $t(\vec K_2, W^{(\lambda)}) = \int_{[0, 1]^2}W^{(\lambda)}(x, y)\;dxdy = 1/16.$ If $B$ does not have a homomorphism to $\vec K_2$, we cannot have a homomorphic copy of $B$ in $W^{(\lambda)}$ where the images of the vertices of $B$, $x_1, \ldots, x_{v(B)}$ lie in $[0, 1/2].$ In other words, when $\lambda = 0$, we have $t(B, W^{(0)}) = 0$. On the other hand, when $\lambda = 1$, we see that the integral is nonzero only when all of $x_i$ lie in $[1/2, 1]$, and $W^{(1)}(x_i, x_j) = 1/4$ in this case. Hence 
\[t(B, W^{(1)}) = \left(\frac{1}{2}\right)^{v(B)}\left(\frac{1}{4}\right)^{e(B)} \ge \left(\frac{1}{16}\right)^{e(B)}.\]
In the last inequality, we note that $B$ has no isolated vertices, so $v(B) \leq 2e(B)$. Since $t(B, W^{(0)}) \leq \left(\frac{1}{16}\right)^{e(B)} \le t(B, W^{(1)})$, by continuity of $t(B, W^{(\lambda)})$ there exists $\lambda_0 \in [0, 1]$ such that
\[t(B, W^{(\lambda_0)}) = \left(\frac{1}{16}\right)^{e(B)} = t\left(\vec K_2, W^{(\lambda_0)}\right)^{e(B)}.\]
However, $W^{(\lambda_0)}$ is not the constant function up to a measure zero set, so $B$ is not forcing.
\end{proof}

Finally, we obtain an analogous result to~\cref{t:disidorenko} for oriented graphs which have the directed forcing property. We show that if $B$ has a homomorphism to an edge, then $\overline{B}$ has the asymmetric forcing property if and only if $B$ has the directed forcing property.

\begin{proof}[Proof of~\cref{t:diforcing}]
First, we show that if $B$ has directed forcing property, then $\overline{B}$ has the asymmetric forcing property.
Suppose that $\{H_n\}_{n\geq 1}$ is a sequence of bipartite graphs such that 
\[\lim_{n\to\infty} t_{\bip}(\ol B, H_n) = p^{e(B)}.\]
For each $n\geq 1$, we set $W_n = \frac12 W_{H_n}$. Then,
\[t_{\bip}(\ol B, W_n) = 2^{-e(B)}t_{\bip}(\ol B, W_{H_n}) = 2^{-e(B)}t_{\bip}(\ol B, H_n).\]
Because $W_n$ satisfies that $W_n(x, y) + W_n(y, x) \leq 1/2 + 1/2 = 1$ for all $x, y\in [0, 1]$, by~\cref{thm:graphon-lim}, there exists a directed graph $G_n$ such that $\delta_\square(W_{G_n}, W_n) \leq 1/n$.
Also note that by~\cref{lem:switching}, $t_{\bip}(\ol B, W_n) = t(B, W_n)$. Hence, 
\[t(B, W_n) = t_{\bip}(\ol B, W_n) = 2^{-e(B)}t_{\bip}(\ol B, H_n).\]
By~\cref{lem:counting}, we have
\[\begin{split}\left|t(B, G_n) - (p/2)^{e(B)}\right| & \leq \left|t(B, W_{G_n}) - t(B, W_n)\right| + \left|2^{-e(B)}t_{\bip}(\ol B, H_n) - 2^{-e(B)} p^{e(B)}\right| \\
& \leq \frac{e(B)}{n} + 2^{-e(B)}\left|t_{\bip}(\ol B, H_n) - p^{e(B)}\right|.\end{split}\]
As $n$ tends to infinity, both terms on the right hand side tend to zero. Therefore we have
\[\lim_{n\to\infty}t(B, G_n) - (p/2)^{e(B)} = 0.\]
By the definition of directed forcing property, we conclude that the sequence of directed graphs $\{G_n\}_{n\geq 1}$ is $(p/2)$-quasirandom. By~\cref{prop:graphon-quasi}, we know that
\[\lim_{n\to\infty}\|W_{G_n} - p/2\|_{\square} = 0.\]
From our construction of $G_n$, we know that $\delta_{\square}(W_{H_n}/2, W_{G_n}) \leq 1/n$. By the triangle inequality for the cut-distance, we have $\delta_{\square}(W_{H_n}/2, p/2) \leq 1/n + \|W_{G_n}-p/2\|_\square$. Note that the constant function is invariant under measure preserving maps, so we can equivalently rewrite this as
\[\|W_{H_n} - p\|_\square \leq 2/n + 2\|W_{G_n}-p/2\|_\square.\]
As $n$ tends to infinity, both terms on the right hand side tend to zero. Therefore, we conclude that
\[\lim_{n\to\infty} \|W_{H_n} - p\| = 0,\]
which establishes $p$-quasirandomness via~\cref{prop:graphon-quasi}.

The other direction is very similar. Suppose that $\ol B$ has the asymmetric forcing property. We aim to show that $B$ has the directed forcing property. Suppose that $\{G_n\}_{n\geq 1}$ is a sequence of directed graphs satisfying that $t(B, G_n) = p^{e(B)}$. Then for $W_{G_n}$~\cref{lem:switching} we have $t_{\bip}(\ol B, W_{G_n} = t(B, W_{G_n}) = t(B, G_n)$.
Note that $W_{G_n}$ is a measurable function taking values in $[0, 1]$, so by~\cref{thm:graphon-lim} there exists a bipartite graph $H_n$ satisfying that $\delta_\square(W_{H_n}, W_{G_n}) \le 1/n$. Using a similar argument as above, we have by~\cref{lem:counting} that
\[\left|t_{\bip}(\ol B, H_n) - p^{e(B)}\right| \le \frac{e(B)}{n} + \left|t(B, G_n) - p^{e(B)}\right|.\]
Letting $B$ tend to infinity, we have $\lim_{n\to\infty}t_{\bip}(\ol B, H_n) = p^{e(B)}.$
Hence $\{H_n\}_{n\geq 1}$ is $p$-quasirandom by the asymmetric forcing property of $\ol B$. This shows that $\lim_{n\to\infty}\|W_{H_n} - p\|_{\square} = 0$. Also by the same argument as before, we have the triangle inequality $\|W_{G_n} - p\|_\square \leq 1/n + \|W_{H_n} - p\|_{\square}$ and thus $\lim_{n\to\infty}\|W_{G_n} - p\|_\square = 0$. By~\cref{prop:graphon-quasi}, we conclude that $\{G_n\}_{n\geq 1}$ is $p$-quasirandom.
\end{proof}

\section{Concluding remarks}\label{s:conc}

The above directed Sidorenko and forcing properties studied here are very natural notions, but not the only Sidorenko-type properties one might investigate in a directed graph.
For example, consider the following, slightly different directed Sidorenko-style property that is natural to define for an oriented graph $B$. We say $B$ satisfies the {\it second directed Sidorenko property} if
\begin{equation}\label{e:sid2}
t_B(G) \ge \frac{1}{2^{e(B)}}  t_{\overline{B}}(\overline{G}) 
\end{equation}
holds for all oriented graphs $G$. 
Equation~\eqref{e:sid2} implies~\eqref{e:sid1} when the underlying graph $\overline{B}$ satisfies Sidorenko's conjecture. This fact suggests that the second directed Sidorenko property captures a natural ``orientation-focused'' analogue of the Sidorenko property. Griffiths~\cite{GRI13}, studied such an orientation-focused analogue of quasirandomness in general oriented graphs.

We can further focus on which patterns of \textit{directions} are systematically overrepresented or underrepresented by studying counts of directed graphs $B$ in \textit{tournaments}, oriented complete graphs. The behavior of subgraph counts in tournaments are quite different than in general directed graphs. For example, there are a small number of special subgraphs, called \textit{impartial oriented graphs} that have the surprising property that for an impartial $B$, the number of copies of $B$ in any tournament $T$ on $n$ vertices is the same, i.e. the number of copies of $B$ in a tournament $T$ is $f(n)$, where $n = v(T)$, and is otherwise independent of the orientation of $T$. Some of the simplest examples of impartial oriented graphs are a vertex, single edge, and the oriented graph $B$ with vertex set $\{a,b,c,d\}$ and edge set $\{(a, b), (c, d), (a, c)\}$ (see~\cite{ZHAO19}).
 
Further, there are directed graphs (e.g. directed paths~\cite{SSZ20} and directed cycles of length not a multiple of $4$~\cite{GKLV20}), that are \textit{tournament anti-Sidorenko}, i.e. systematically underrepresented in all tournaments. This phenomenon does not appear in general directed graphs. The tournament Sidorenko and tournament anti-Sidorenko directed properties are further studied in~\cite{FHMZ22}.


\bibliographystyle{alpha}
\bibliography{refs.bib}

\newcommand{\etalchar}[1]{$^{#1}$}
\begin{thebibliography}{HKK{\etalchar{+}}19}

\bibitem[AGH11]{AM11}
Omid Amini, Simon Griffiths, and Florian Huc.
\newblock Subgraphs of weakly quasi-random oriented graphs.
\newblock {\em SIAM J. Discrete Math.}, 25(1):234--259, 2011.

\bibitem[BCL{\etalchar{+}}08]{BCL08}
C.~Borgs, J.~T. Chayes, L.~Lov\'{a}sz, V.~T. S\'{o}s, and K.~Vesztergombi.
\newblock Convergent sequences of dense graphs. {I}. {S}ubgraph frequencies,
  metric properties and testing.
\newblock {\em Adv. Math.}, 219(6):1801--1851, 2008.

\bibitem[BLSS21]{BU19}
Matija Buci\'{c}, Eoin Long, Asaf Shapira, and Benny Sudakov.
\newblock Tournament quasirandomness from local counting.
\newblock {\em Combinatorica}, 41(2):175--208, 2021.

\bibitem[CFS10]{CON10}
David Conlon, Jacob Fox, and Benny Sudakov.
\newblock An approximate version of {S}idorenko's conjecture.
\newblock {\em Geom. Funct. Anal.}, 20(6):1354--1366, 2010.

\bibitem[CG91]{FAN91}
F.~R.~K. Chung and R.~L. Graham.
\newblock Quasi-random tournaments.
\newblock {\em J. Graph Theory}, 15(2):173--198, 1991.

\bibitem[CGW89]{FAN89}
F.~R.~K. Chung, R.~L. Graham, and R.~M. Wilson.
\newblock Quasi-random graphs.
\newblock {\em Combinatorica}, 9(4):345--362, 1989.

\bibitem[CKLL18]{CON18}
David Conlon, Jeong~Han Kim, Choongbum Lee, and Joonkyung Lee.
\newblock Some advances on {S}idorenko's conjecture.
\newblock {\em J. Lond. Math. Soc. (2)}, 98(3):593--608, 2018.

\bibitem[CPS19]{COR15}
Leonardo~N. Coregliano, Roberto~F. Parente, and Cristiane~M. Sato.
\newblock On the maximum density of fixed strongly connected subtournaments.
\newblock {\em Electron. J. Combin.}, 26(1):Paper No. 1.44, 48, 2019.

\bibitem[CR17]{COR17}
Leonardo~Nagami Coregliano and Alexander~A. Razborov.
\newblock On the density of transitive tournaments.
\newblock {\em J. Graph Theory}, 85(1):12--21, 2017.

\bibitem[FHMZ22]{FHMZ22}
Jacob Fox, Zoe Himwich, Nitya Mani, and Yunkun Zhou.
\newblock Tournament anti-{S}idorenko digraphs.
\newblock {\em in preparation}, 2022.

\bibitem[GKLV20]{GKLV20}
Andrzej Grzesik, Daniel Kral, Laszlo~Miklos Lovasz, and Jan Volec.
\newblock Cycles of a given length in tournaments.
\newblock {\em arXiv preprint arXiv:2008.06577}, 2020.

\bibitem[Gri13]{GRI13}
Simon Griffiths.
\newblock Quasi-random oriented graphs.
\newblock {\em J. Graph Theory}, 74(2):198--209, 2013.

\bibitem[HKK{\etalchar{+}}19]{HA19}
Robert Hancock, Adam Kabela, Daniel Kral, Taisa Martins, Roberto Parente, Fiona
  Skerman, and Jan Volec.
\newblock No additional tournaments are quasirandom-forcing.
\newblock {\em arXiv preprint arXiv:1912.04243}, 2019.

\bibitem[KLL16]{KIM16}
Jeong~Han Kim, Choongbum Lee, and Joonkyung Lee.
\newblock Two approaches to {S}idorenko's conjecture.
\newblock {\em Trans. Amer. Math. Soc.}, 368(7):5057--5074, 2016.

\bibitem[Lov12]{LOV12}
L\'{a}szl\'{o} Lov\'{a}sz.
\newblock {\em Large networks and graph limits}, volume~60 of {\em American
  Mathematical Society Colloquium Publications}.
\newblock American Mathematical Society, Providence, RI, 2012.

\bibitem[LS11]{LI17}
J.L.~Xiang Li and Bal\'azs Szegedy.
\newblock On the logarithimic calculus and {S}idorenko's conjecture.
\newblock {\em arXiv preprint arXiv:1107.1153}, 2011.

\bibitem[Sid93]{SID93}
Alexander Sidorenko.
\newblock A correlation inequality for bipartite graphs.
\newblock {\em Graphs Combin.}, 9(2):201--204, 1993.

\bibitem[Sim84]{SIM82}
Mikl{\'o}s Simonovits.
\newblock Extremal graph problems, degenerate extremal problems, and
  supersaturated graphs.
\newblock {\em Progress in graph theory (Waterloo, Ont., 1982)}, pages
  419--437, 1984.

\bibitem[SSZ20]{SSZ20}
Ashwin Sah, Mehtaab Sawhney, and Yufei Zhao.
\newblock Paths of given length in tournaments.
\newblock {\em arXiv preprint arXiv:2012.00262}, 2020.

\bibitem[ST04]{SKO04}
Jozef Skokan and Lubos Thoma.
\newblock Bipartite subgraphs and quasi-randomness.
\newblock {\em Graphs Combin.}, 20(2):255--262, 2004.

\bibitem[Sze14]{SZE}
Bal\'azs Szegedy.
\newblock An information theoretic approach to {S}idorenko's conjecture.
\newblock {\em arXiv preprint arXiv:1406.6738}, 2014.

\bibitem[Tho87]{TH87}
Andrew Thomason.
\newblock Pseudo-random graphs.
\newblock {\em Ann. Discrete Math. (33)}, 144:307--331, 1987.

\bibitem[ZZ20]{ZHAO19}
Yufei Zhao and Yunkun Zhou.
\newblock Impartial digraphs.
\newblock {\em Combinatorica}, 40(6):875--896, 2020.

\end{thebibliography}

\end{document}